\magnification=\magstep1
\def\ds{\baselineskip 20pt plus 2pt}
\def\ss{\baselineskip 10pt plus 1pt}
\ss
\input amssym
\def\N{\Bbb N}
\voffset=1in
\def\ans{\vrule height.1pt width80pt depth0pt}
\def\nuc{\operatorname{nuc}}
\pageno=1
\footline={\ifnum\pageno=1\hfill\else\hss\tenrm\folio\hss\fi}
\font\sl=cmbsy10
\def\slN{\hbox{$N \textfont1=\sl$}}

\font\bigbold=cmbx10 scaled \magstep2
\centerline {\bigbold Nuclear Operators on Spaces}
\centerline {\bigbold of Continuous Vector-Valued Functions}
\vskip.75truein
\centerline {by}
\vskip.75truein
\centerline {\bf Paulette Saab$^*$ and Brenda Smith}
\vskip 1truein
{\narrower\smallskip\noindent {\bf Abstract}\ \ Let $\Omega$ be a
compact Hausdorff space, let $E$ be a Banach space, and let $C(\Omega,
E)$ stand for the Banach space of all $E$-valued continuous functions on
$\Omega$ under supnorm.  In this paper we study when nuclear operators
on $C(\Omega, E)$ spaces can be completely characterized in terms of
properties of their representing vector measures.  We also show that if
$F$ is a Banach space and if $T:\ C(\Omega, E)\rightarrow F$ is a
nuclear operator, then $T$ induces a bounded linear operator $T^\#$ from
the space $C(\Omega)$ of scalar valued continuous functions on $\Omega$
into $\slN(E,F)$ the space of nuclear operators from $E$ to $F$, in this
case we show that $E^*$ has the  Radon-Nikodym property if and
only if $T^\#$ is nuclear whenever $T$ is nuclear.
\vskip2truein
\noindent AMS(MOS) Subject Classification (1980)\hfill\break
Primary 46E40, 46G10, 47B10, Secondary 28B05, 28B20

\noindent \ans

\noindent $^*$ Research supported in part by an NSF Grant DMS-87500750}
\vfill\eject
\voffset=-.3in
\ds
\
\vskip .4truein
\
Let $\Omega$ be a compact Hausdorff space, let $E$ be a Banach space,
and let $C(\Omega, E)$ stand for the Banach space of continuous
$E$-valued functions on $\Omega$ under supnorm.  It is well known [3,
p.182] that if $F$ is a Banach space, then any bounded linear operator
$T:\ C(\Omega, E)\longrightarrow F$ has a finitely additive vector
measure $G$ defined on the $\sigma$-field of Borel subsets of $\Omega$
with values in the space {\it \$}$(E,F^{**})$ of bounded linear
operators from $E$ to the second dual $F^{**}$ of $F$.  The measure $G$
is said to represent $T$.  The purpose of this note is to study the
interplay between certain properties of the operator $T$ and properties
of the representing measure $G$.  Precisely, one of our goals is to
study when one can characterize nuclear operators in terms of their
representing measures.  This is of course motivated by a well known
Theorem of L. Schwartz [5] (see also [3, p.173]) concerning nuclear
operators on spaces $C(\Omega)$ of continuous scalar-valued functions.
The study of nuclear operators on spaces $C(\Omega, E)$ of continuous
vector-valued functions was initiated in [1] where the author extended
Schwartz's result in case $E^*$ has the  Radon-Nikodym
property.  In this paper, we will show that the condition on $E^*$ to
have the Radon-Nikodym property  is necessary to have a
Schwartz's type theorem. This leads to a new characterization of dual
spaces $E^*$ with the Radon-Nikodym property.
In [2] it was shown that if $T:\
C(\Omega,E)\rightarrow F$ is nuclear then its representing measure $G$
takes its values in the space $\slN(E,F)$ of nuclear operators
from $E$ to $F$.  One of the results of this paper is that if $T:\ \
C(\Omega, E)\rightarrow F$ is nuclear then its representing measure $G$
is countably additive and
of bounded variation as a vector measure taking its values in
$\slN(E,F)$ equipped with the nuclear norm.
Finally, we show by easy examples that the above mentioned conditions on
the representing measure $G$ do not characterize nuclear operators on
$C(\Omega,E)$ spaces, and we also look at cases where nuclear operators
are indeed characterized by the above two conditions.  For all undefined
notions and terminologies we refer the reader to [3].
\bigskip
\noindent {\bf 0--Preliminaries}\ \ If $X$ and $Y$ are Banach spaces,
then {\it \$}$(X,Y)$ will stand for the space of bounded linear
operators from $X$ to $Y$.  An element $T$ in {\it \$}$(X,Y)$ is said to
be a {\bf nuclear operator} if there exist sequences $(x^*_n)$ in $X^*$
and $(y_n)$ in $Y$ such that for each $x$ in $X$
$$
T(x)=\sum\limits^\infty_{n=1}x_n^*(x)y_n,
$$
and
$$
\sum\limits^\infty_{n=1}\parallel x_n^*\parallel\ \parallel
y_n\parallel<\infty.
$$

We say that $\sum\limits_n x_n^*\otimes y_n$ represents the nuclear
operator $T$.  The {\bf nuclear norm} of a nuclear operator $T:\
X\rightarrow Y$ is defined by:
$$
\parallel T\parallel_{\nuc}=\inf\left\{\sum\limits_n\parallel
x_n^*\parallel\ \parallel y_n\parallel\right\}
$$
where the infimum is taken over all sequences $(x_n^*)$ and $(y_n)$ such
that $T(x)=\sum\limits^\infty_{n=1}x_n^*(x)y_n$ holds for all $x$ in
$X$.  The nuclear operators from $X$ to $Y$ form a normed linear space
under the nuclear norm [3, p.170], which we shall denote by $\slN(X,Y)$.

If
$\Omega$ is a compact Hausdorff space and $E$ is a Banach space, then
$C(\Omega, E)$ will stand for the space of continuous $E$-valued
functions defined on $\Omega$ under supnorm.  If $E=\Bbb R$ or $\Bbb C$
we will simply write $C(\Omega)$.  The space $M(\Omega, E^*)$ will stand
for the space of all regular $E^*$-valued vector measures $\mu$ defined
on the $\sigma$-field $\sum$ of Borel subsets of $\Omega$ that are of
bounded variation.  We shall use the fact (see [3, p.182]) that
$M(\Omega, E^*)$ is a Banach space under the variation norm $\parallel
\mu \parallel=|\mu|(\Omega)$, and that $M(\Omega, E^*)$ is isometrically
isomorphic to the dual space $C(\Omega, E)^*$.  When $E=\Bbb R$ or $\Bbb
C$ we will simply write $M(\Omega)$.  If $\mu \in M(\Omega, E^*)$, then
for each $e\in E$ we will denote by $\langle e,\mu\rangle$ the element of
$M(\Omega)$ such that for each $f\in C(\Omega)$,
$$
\int f d\langle e,\mu\rangle =\mu (f\otimes e)
$$
where $f\otimes e$ is the element in $C(\Omega, E)$ such that $f\otimes
e(\omega)=f(\omega) e$ for each $\omega\in \Omega$.

If $\nu \in M(\Omega)$ and $x^*\in
E^*$, we denote by $\nu \otimes x^*$ the element of $M(\Omega, E^*)$
which to each Borel subset $B$ of $\Omega$ associates the element $\nu
(B)x^*$ of $E^*$. lf $E$ and $F$ are Banach spaces and
$\Omega$ is a compact Hausdorff space, then we will denote by $G$ the
finitely additive {\it \$}$(E,F^{**})$-valued measure representing the
operator $T$.  Recall that if $B$ is a Borel subset of $\Omega$
then
$$
G(B)e=T^{**}(\phi_{B,e})
$$
for all $e\in E$, where $\phi_{B,e}$ is the element of $M(\Omega, E^*)^*$
such for each $\lambda \in M(\Omega, E^*)$
$$
\phi_{B,e}(\lambda)=\lambda(B)(e)
$$
and $T^{**}$ is the second adjoint of $T$.

Finally we recall that a Banach space $X$ has the {\bf Radon-Nikodym
property} (RNP) if for every finite measure space $(S, \Sigma, \mu)$ and
every vector measure $m:\ \Sigma \rightarrow X$ of bounded variation that
is absolutely continuous with respect to $\mu$ there exists a strongly
measurable Bochner integrable function $g:\ S\rightarrow X$ such that
$$
m(A)=\int_Afd\mu
$$
for each $A\in \Sigma$.
\bigskip
\noindent {\bf 1--Some Properties of the Measure representing a Nuclear
operator}
\bigskip
Throughout we let $\Omega$ be a compact Hausdorff space with $\sum$ its
$\sigma$-field of Borel subsets and we let $E$ and $F$ be Banach spaces.
In what follows we shall look at some of the properties that a nuclear
operator on $C(\Omega, E)$ induces on its representing measure $G$.  In
[2] it was shown that if $T:\ C(\Omega, E)\rightarrow F$ is nuclear then
$G:\ \sum \rightarrow \slN(E,F)$.  In the next proposition we shall show
that $G$ enjoys a stronger property.
\vfill\eject
\bigskip
\noindent {\bf Proposition 1}\ \ If $T:\ C(\Omega, E)\rightarrow F$ is a
nuclear operator with representing measure $G$, then:

\itemitem {(i)} for each Borel subset $B$ of $\Omega,\ G(B):\
E\rightarrow F$ is a nuclear operator, and

\itemitem {(ii)} the measure $G$ is countably additive and is of bounded
variation as a
vector measure taking its values in $\slN(E,F)$ under the nuclear norm.

\medskip
\noindent {\bf Proof:}\ \ If $T:\ C(\Omega, E)\rightarrow F$ is nuclear,
then there are sequences $(\mu_n)$ in $M(\Omega, E^*)$ and $(y_n)$ in
$F$ such that for all $f\in C(\Omega, E)$
$$
Tf=\sum\limits^\infty_{n=1}\mu_n(f)y_n
$$
and
$$
\sum\limits^\infty_{n=1} \parallel \mu_n\parallel \, \parallel
y_n\parallel \,
<\infty.
$$
In particular the operator $T^{**}$ is nuclear with
$$
T^{**}(\phi)=\sum\limits^{\infty}_{n=1}\phi (\mu_n)y_n
$$
for all $\phi \in C(\Omega, E)^{**}.$  This implies that for any Borel
subset $B$ of $\Omega$ and for all $e\in E$
$$
\aligned
G(B)e&=T^{**}(\phi_{B,e})\\
&=\sum\limits^\infty_{n=1}\mu_n (B)(e)y_n\endaligned
$$
Since $\left(\mu_n(B)\right)$ is a sequence in $E^*$ and
$\sum\limits^\infty_{n=1}\parallel \mu_n (B)\parallel\ \parallel
y_n\parallel \leq \sum\limits^\infty_{n=1}\parallel \mu_n\parallel\
\parallel y_n\parallel <\infty$, this will quickly show that for each
Borel subset $B$ of $\Omega$, the operator $G(B):\ E\rightarrow F$ is
nuclear.  To prove (ii) note that since $G:\
\sum\rightarrow \slN (E,F)$, let $|G|_{\nuc}$ denote
the extended non negative function whose value on a set $B$ in $\sum$ is
given by
$$
|G|_{\nuc}(B)=\sup\limits_\pi \sum\limits_{B_i\in \pi} \parallel
G(B_i)\parallel_{\nuc}
$$
where the supremum is taken over all finite partitions $\pi$ of $B$.  We
first show that \hfill\break $|G|_{\nuc}(\Omega)<\infty.$

For this, note that if $\pi=\{B_i\}$ is a finite partition of $\Omega$,
then
$$
\aligned
\sum\limits_{B_i\in \pi}\parallel G(B_i)\parallel_{\nuc}&\leq
\sum\limits_{B_i\in \pi} \sum\limits^\infty_{n=1} \parallel \mu_n
(B_i)\parallel\ \parallel y_n\parallel\\
&\leq \sum\limits_{B_i\in \pi} \sum\limits^\infty_{n=1} |\mu_n|(B_i)
\parallel y_n\parallel\\
&\leq \sum\limits^\infty_{n=1}\sum\limits_{B_i\in \pi} |\mu_n|
(B_i)\parallel y_n\parallel\\
&\leq \sum\limits^\infty_{n=1} |\mu_n|(\Omega)\parallel y_n\parallel\\
&=\sum\limits^\infty_{n=1}\parallel \mu_n\parallel\ \parallel
y_n\parallel\quad <\infty\endaligned
$$
This of course shows that $|G|_{\nuc}(\Omega)<\infty$.  To
complete the proof of (ii) we need to show that $G$ is countably
additive.  First note that since $|G|_{\nuc}(\Omega)<\infty$, it follows
from [3, p.7] that $G$ is strongly additive, that is if $(B_i)$ is a
sequence of pairwise disjoint Borel subsets of $\Omega$ we have that the
series $\sum\limits^\infty_{i=1}G(B_i)$ converges in $\slN(E,F)$.  To
complete the proof we need to check that the series
$\sum\limits^\infty_{i=1}G(B_i)$ converges to
$G(\dsize\bigcup\limits_{i\geq 1}B_i)$ in $\slN (E,F)$.  To this end,
consider a
series $\sum\limits_{n=1}^\infty \mu_n \otimes y_n$ representing the
operator
$T$
such that $$\sum\limits^\infty_{n=1}||\mu_n|| \, ||y_n||<\infty .$$
Without
loss of generality we may and shall assume that $||y_n||\leq 1$ for all
$n\geq 1$.  Let $\epsilon >0$ and pick $N\in \N$ such that
$$
\sum\limits^\infty_{\nu=N +1}||\mu_n||<\epsilon/2
$$

Since $\mu_1, \mu_2,\ldots$, and $\mu_N\in M(\Omega, E^*)$, there
exists $K\in \N$ such that
$$
|\mu_j|\left(\bigcup_{i\geq K}B_i\right)<\dfrac \epsilon {2^{j+1}} \text
{ for all }j=1,\ldots, N
$$
This implies that
$$
\aligned
||G(\bigcup\limits_{i\geq
1}B_i)-\sum\limits^{K-1}_{i=1}G(B_i)||_{\nuc}&=\\
||G\left(\bigcup\limits_{i\geq K}B_i\right)||_{\nuc}&\leq\\
\sum\limits^\infty_{n=1}||\mu_n\left(\bigcup_{i\geq
K}B_i\right)||&\leq\endaligned
$$

$$
\sum\limits^N_{n=1}|\mu_n|\left(\bigcup\limits_{i\geq K}B_i\right)+
\sum\limits^\infty_{n=N+1}|\mu_n|\left(\bigcup\limits_{i\geq
K}B_i\right)<\epsilon
$$
This shows that $G$ is countably additive as a vector measure taking its
values in $\slN (E,F)$.

The first question that arises at this stage, is when do properties (i)
and (ii) above characterize nuclear operators?  The next proposition
shows that one necessary condition is that $F$ should have the
Radon-Nikodym property.
\medskip
\noindent {\bf Proposition 2}\ \ If $F$ fails RNP, then for any Banach
space $E$ there is a non nuclear operator $T:\ C([0,1], E)\rightarrow F$
whose representing measure takes its values in $\slN(E,F)$ and is of
bounded variation as a countably additive vector measure taking its
values in
$\slN(E,F)$.
\medskip
\noindent {\bf Proof:}\ \ If $F$ lacks RNP then by [3, p.175] there is
an $F$-valued countably additive vector measure $m$ on the Borel subsets
of [0,1] such that $m$ is of bounded variation, $m$ is absolutely
continuous with respect to Lebesgue measure but $m$ admits no Bochner
integrable derivative with respect to its variation $|m|$.  If we
define $T':\ C[0,1]\rightarrow F$ by
$$
T'(f)=\int_{[0,1]}f\,dm
$$
then $T'$ is not a nuclear operator (see [5] or [3, p.173]).  Now fix
$e\neq 0$ in $E$ then choose $e^*$ in $E^*$ with $e^*(e)=1$ and define
$T:\ C([0,1],E)\rightarrow F$ by
$$
T(\varphi)=T'(e^*\circ \varphi)
$$
for each $\varphi$ in $C([0,1], E)$.  In particular, for each $f$ in
$C[0,1]$ and $x$ in $E$ we have
$$
T(f\otimes x)=e^*(x)T'(f).
$$
It is clear that the measure representing the operator $T$ is the
measure $G=m\otimes e^*$ which to every Borel subset $B$ of [0,1]
associates the one-rank operator such that $G(B)x=e^*(x)m(B)$ for each
$x\in E$.  On other hands for any finite partition $\pi=\{B_i\}$ of
[0,1] we have
$$
\aligned
\sum\limits_{B_i\in \pi}\parallel G(B_i)\parallel_{\nuc}&\leq
\sum\limits_{B_i\in \pi}\parallel e^*\parallel \parallel m
(B_i)\parallel\\
&\leq \parallel e^*\parallel\ |m|([0,1])<\infty.\endaligned
$$
so $$|G|_{\nuc}([0,1])<\infty.$$
If $T$ were a nuclear operator then there would exist $(\mu_n)$ in
$M([0,1], E^*)$ and $(y_n)$ in $F$ such that for each $\phi$ in
$C([0,1], E)$
$$
T(\phi)=\sum\limits^\infty_{n=1}\mu_n(\phi) y_n
$$
In particular for each $f\in C[0,1]$
$$
\aligned
T(f\otimes e)&=T'(f)\\
&=\sum\limits^\infty_{n=1}\mu_n (f\otimes e) y_n\\
&=\sum\limits^\infty_{n=1} \int_{[0,1]}f d\langle e, \mu_n\rangle
y_n\endaligned
$$
This of course shows that $T'$ is represented by $\sum\limits_n \langle
e, \mu_n\rangle \otimes y_n$; moreover
$$
\sum\limits^\infty_{n=1}\parallel
\langle e,
\mu_n \rangle \parallel \parallel y_n\parallel \leq \parallel e\parallel
\sum\limits_n \parallel \mu_n\parallel \parallel y_n\parallel<\infty
$$
which implies that $T'$ is nuclear. This contradiction  finishes
the proof.
\bigskip
This brings us to ask the following question.  Let $E$ and $F$ be Banach
spaces such that $F$ has the Radon-Nikodym property.  Let $\Omega $ be a
compact Hausdorff space and let $T:\ C(\Omega, E)
\rightarrow F$ be a bounded linear operator satisfying conditions (i)
and (ii) of Proposition 2.  Is $T$ nuclear?
Recently, the first named author has
given a positive answer to the above question when $F$ is complemented
in its bidual $F^{**}$ (see [4]).
\bigskip
\noindent {\bf 2--The Radon-Nikodym Property and Nuclear Operators}

Throughout this section $\Omega$ is a compact Hausdorff space and $E$
and $F$ are Banach spaces.  Every bounded linear operator $T:\ C(\Omega,
E)\rightarrow F$ induces a bounded linear operator $T^\#:\
C(\Omega)\longrightarrow${\it \$}$(E,F)$ where for each $f$ in
$C(\Omega)$
$$
T^\#(f)(e)=T(f\otimes e)
$$
for all $e\in E$.

In this section we shall look at the interplay of the two operators $T$
and $T^\#$.  The next result shows that when $T$ is nuclear, then the
range of $T^\#$ is  in $\slN(E,F)$.
\bigskip
\noindent {\bf Theorem 3}\ \ If $T:\ C(\Omega, E)\rightarrow F$ is nuclear,
then $T^\#$ takes its values in $\slN(E,F).$
\medskip
\noindent {\bf Proof}\ \ Let $(\mu_n)$ in $M(\Omega, E^*)$ and $(y_n)$ in
$F$ be such that $\sum\limits^\infty_{n=1}\parallel \mu_n\parallel
\parallel y_n\parallel <\infty$ and for each $\varphi$ in $C(\Omega, E)$
$$
T(\varphi)=\sum\limits^\infty_{n=1}\mu_n(\varphi)y_n.
$$
For each $n\geq 1$ define $\mu^\#_n:\ C(\Omega)\rightarrow E^*$ by
$$\mu_n^\#(f)e=\mu_n(f\otimes e)
$$
for each $f\in C(\Omega)$ and $e\in E$.
Since for each $f\in C(\Omega)$ and $e\in E$
$$
\aligned
T^\#(f)(e)&=T(f\otimes e)\\
&=\sum\limits^\infty_{n=1}\mu_n(f\otimes e)y_n\\
&=\sum\limits^\infty_{n=1} \mu^\#_n(f)(e)y_n\endaligned
$$
it follows that for each $f$ in $C(\Omega)\ T^\#(f)$ can be represented
by
the series $\sum\limits^\infty_{n=1}\mu^\#_n(f)\otimes y_n$.  Moreover
since $\sum\limits^\infty_{n=1}\parallel \mu_n^\#(f)\parallel \parallel
y_n\parallel \leq \parallel f\parallel \sum\limits^\infty_{n=1}\parallel
\mu_n\parallel \parallel f_n\parallel <\infty$, it follows that
$T^\#(f)$ is a nuclear operator from $E$ to $F$ for each $f\in
C(\Omega)$.

The next result illustrates one key relationship between $T$ and $T^\#$.
\bigskip
\noindent {\bf Theorem 4}\ \ The operator $T\ C(\Omega, E)\rightarrow
F$ is nuclear whenever $T^\#:\ C(\Omega)\rightarrow \slN(E,F)$ is
nuclear.
\bigskip
\noindent {\bf Proof:}\ \ Assume $T^\#:\ \ C(\Omega)\rightarrow
\slN(E,F)$ is nuclear. Then there exist sequences $(\nu_n)$ in
$C(\Omega)^*$
and $(N_n)$ in $\slN(E,F)$ such that for each $f\in C(\Omega)$
$$
T^\#(f)=\sum\limits^\infty_{n=1}\nu_n(f)N_n
$$
and
$$
\sum\limits^\infty_{n=1}\parallel \nu_n\parallel \parallel
N_n\parallel_{\nuc}<\infty.
$$
Without loss of generality we may and do assume that $\parallel
\nu_n\parallel \leq 1$ for all $n\geq 1.$
Similarly, since
each $N_n$ is a nuclear operator, for each $n\geq 1$, there are
sequences $(e^*_{n,m})$ in $E^*$ and $(y_{n,m})$ in $F$ such that for
all $e\in E$
$$
N_n(e)=\sum\limits^\infty_{m=1}e^*_{n,m}(e) y_{n,m}
$$
and
$$
\sum\limits^\infty_{m=1}\parallel e^*_{n,m}\parallel \parallel
y_{n,m}\parallel \leq \parallel N_n\parallel_{\nuc}+1/2^n.
$$
Hence
$$
\sum\limits^\infty_{n=1}\sum\limits^\infty_{m=1}\parallel \nu_n\parallel
\parallel e^*_{n,m}\parallel \parallel y_{n,m}\parallel \leq
\sum\limits^\infty_{n=1}\parallel \nu_n\parallel \parallel
N_n\parallel_{\nuc}+\sum\limits^\infty_{n=1}\parallel
\nu_n\parallel\big/2^n<\infty \tag {$\dagger\dagger$}
$$
Moreover, for each $f\in C(\Omega)$ and $e\in E$
$$
\aligned
T(f\otimes e)&=T^\#(f)(e)\\
&=\left(\sum\limits^\infty_{n=1}\nu_n(f)N_n\right)(e)\\
&=\sum\limits^\infty_{n=1}\sum\limits^\infty_{m=1}\nu_n(f)e^*_{n,m}(e)y_{n,m}
\\
&=\sum\limits^\infty_{n=1}\sum\limits^\infty_{m=1}\nu_n\otimes
e^*_{n,m}(f\otimes e)y_{n,m}\endaligned
$$
Since the set $\{f\otimes e:\ f\in C(\Omega), e\in E\}$ is total in
$C(\Omega, E)$,  we can assert
that $T$ is represented by the double indexed series
$\sum\limits_n \sum\limits_m(\nu_n\otimes e^*_{n,m})\otimes y_{n,m}$,
where $(\nu_n\otimes e_{n,m}^*)$ is in $M(\Omega, E^*)\simeq C(\Omega,
E)^*$ and $(y_{n,m})$ is in $F$. An appeal to
($\dagger\dagger$)
shows that $T$ is nuclear.
\medskip
In the following we will show that the converse of Theorem 4 does not
always hold.  But if in addition $E^*$ is assumed to have the
Radon-Nikodym property then a close look at [1, Theorem III.4] reveals
that any nuclear operator $T:\ C(\Omega, E)\rightarrow F$ will indeed
induce a nuclear operator $T^{\#}:\ C(\Omega)\rightarrow \slN(E,F)$.
Moreover our next result shows how critical the condition on $E^*$ to
have the Radon-Nikodym property in order that $T$ nuclear implies
$T^{\#}$ nuclear.  As a matter of fact one can characterize dual Banach
spaces with the Radon-Nikodym property as follows:

\noindent {\bf Theorem 5}\ \ Let $E$ and $F$ be Banach spaces.  The
following properties are equivalent:

\itemitem {(i)} The dual space $E^*$ has (RNP);

\itemitem {(ii)} For every compact space $\Omega$, a bounded linear
operator $T:\ C(\Omega, E)\rightarrow F$ is nuclear if and only if
$T^\#:\ C(\Omega)\rightarrow \slN(E,F)$ is nuclear.
\medskip
\noindent {\bf Proof:}\ \ (ii)$\Rightarrow$(i).

If $E^*$ fails RNP, by [3, p.175] there exists $\mu\in M([0,1],E^*)$
such that the operator $\mu^\#:\ C[0,1]\rightarrow E^*$ defined by
$$
\mu^\#(g)e=\mu(g\otimes e)
$$
for all $g\in C[0,1]$ and $e\in E$, is not nuclear.  Choose $y\in F$
such that $y\neq 0$, and define $T:\ C([0,1], E)\rightarrow F$ by
$$
T(\varphi)=\mu(\varphi)y
$$
for all $\varphi\in C([0,1], E)$.  It is clear that $T$ is a rank-one
operator, hence it is nuclear.  The operator $T^\#:\ C[0,1]\rightarrow
\slN(E,F)$ induced by $T$ is clearly the operator such that for each
$g\in C[0,1]$
$$
T^\#(g)=\mu^\#(g)\otimes y.
$$
To see that $T^\#$ is not nuclear, note that for each $y^*\in F^*$, we
can define the operator $T^\#_{y^*}:\ C[0,1]\rightarrow E^*$ by
$$
T^\#_{y^*}(g)=(T^\# g)^*(y^*)
$$
for each $g\in C[0,1]$.

If $T^\#$ were a nuclear operator then $T^\#_{y^*}$ would also be a
nuclear operator for each $y^*\in F^*$.  This of course follows from the
fact that $T^\#_{y^*}$ is the composition of $T^\#$ and the bounded
linear operator from $\slN(E,F)$ to $E^*$ which to an element $N$ in
$\slN(E,F)$ associates the element $N^*(y^*)$ in $E^*$.  But for each
$g$ in $C[0,1]$
$$
T^\#_{y^*}(g)=y^*(y)\mu^\#(g).
$$
By the Hahn-Banach Theorem, choose $y^*$ in $F^*$ such that $y^*(y)=1$,
then for this particular $y^*$ we have
$$
T^\#_{y^*}=\mu^\#.
$$
This contradiction shows that $T^\#$ can not be nuclear.  This proves
(ii)$\Rightarrow $(i).

The proof of (i)$\Rightarrow$ (ii) is implicit in [1, Theorem III.4]. We
shall provide a sketch of a proof for the sake of completeness.  For
this, assume that $T:\ C(\Omega, E)\rightarrow F$ is nuclear and that
$E^*$has the Radon-Nikodym property.
By Proposition 1 we know that the
measure $G$ representing the operator $T$ is countably additive
as a vector measure taking its values in
$\slN(E,F)$ and $|G|_{\nuc}(\Omega)<\infty$.  Here it is easy to
note that it follows from general vector measure techniques [3, p.3]
that
$|G|_{\nuc}$ is countably additive.  The proof that $T^\#$ is
nuclear now follows the proof of the scalar case as given in [3, p.173]
with some minor changes.  For instance, since $E^*$ has the
Radon-Nikodym property, one
can proceed to produce a Bochner $|G|_{\nuc}$-integrable function
$$
H:\ \Omega \longrightarrow \slN (E,F)
$$
such that for each Borel subset $B$ of $\Omega$,
$$
G(B)=\int_BH(\omega) d|G|_{\nuc}(\omega).
$$
and for each $f\in C(\Omega)$ and each $e\in E$
$$
T(f\otimes e)=\int_\Omega f(\omega) H(\omega)(e) d|G|_{\nuc}
(\omega).
$$
Hence for each $f\in C(\Omega)$
$$
T^\#(f)=\int_\Omega f(\omega) H(\omega) d|G|_{\nuc}(\omega)
$$
Another appeal to [3, p.173] shows that $T^\#$ is nuclear.
\vfill\eject
\centerline {\bf Bibliography}
\bigskip
\item {[1]} G. Alexander, {\it Linear Operators on the Spaces of
Vector-Valued Continuous Functions}, Ph. Dissertation, New Mexico State
University, Las Cruces, New Mexico, 1972.

\item {[2]} R. Bilyeu and P. Lewis, {\it Some mapping properties of
representing measures}, Ann. Mat. Pure. Appl. CIX (1976) pp. 273--287.

\item {[3]} J. Diestel, and J.J. Uhl, Jr. {\bf Vector Measures}, Math
Surveys, {\bf 15}, AMS, Providence, RI (1977).

\item {[4]} P. Saab, {\it Integral Operators on Spaces of Continuous
Vector-Valued Functions}, Proc. Amer. Math. Soc. (to appear).

\item {[5]} L. Schwartz, {\it S\'eminaire Schwartz, Expos\'e 13},
Universit\'e de Paris, Facult\'e des Sciences, (1953/54).
\vskip 1truein
\centerline {\bf University of Missouri}
\centerline {\bf Mathematics Department}
\centerline {\bf Columbia, MO  65211}

\vfill\eject\end